\documentclass[runningheads]{llncs}
\usepackage{amssymb,enumerate,tikz,float,array,commath,multirow}
\usetikzlibrary{matrix,arrows,decorations.pathmorphing,automata,positioning}
\usetikzlibrary{decorations.markings}
\usepackage{amsfonts,bbding}
\usepackage{bbm}
\def\natnum{{\mathbb N}}

\def\email#1{Email: {\tt #1}}
\def\proof{\noindent{\bf Proof.}\enspace}

\def\qed{~~\hbox{\hskip 1pt \vrule width 4pt height 8pt depth 1.5pt\hskip 1pt}}
\spnewtheorem{thm}{Theorem}{\bf }{\it }
\spnewtheorem{prop}[thm]{Proposition}{\bf }{\it }
\spnewtheorem{prob}[thm]{Open Problem}{\bf }{\it }
\spnewtheorem{cor}[thm]{Corollary}{\bf }{\it }
\spnewtheorem{lem}[thm]{Lemma}{\bf }{\it }
\spnewtheorem{defn}[thm]{Definition}{\bf }{\rm }
\spnewtheorem{rem}[thm]{Remark}{\bf }{\rm }
\spnewtheorem{exmp}[thm]{Example}{\bf }{\rm }
\spnewtheorem{clm}[thm]{Claim}{\bf }{\it }
\spnewtheorem{nota}[thm]{Notation}{\bf }{\rm }
\spnewtheorem{ppt}[thm]{Property}{\bf }{\rm }
\spnewtheorem{fact}[thm]{Fact}{\bf }{\rm }

\newcommand{\charf}{{\mathbbm{1}}}

\newcommand{\RTD}{\mbox{RTD}}
\newcommand{\TS}{\mbox{TS}}
\newcommand{\TP}{\mbox{TP}}
\newcommand{\TD}{\mbox{TD}}
\newcommand{\XTD}{\mbox{XTD}}

\newcommand{\cL}{{\mathcal{L}}}
\newcommand{\cC}{{\mathcal{C}}}

\newcommand\spn[1]{{\left\langle#1\right\rangle}}

\newcommand{\sm}{\setminus}

\newcommand{\DS}{\mbox{DS}}

\newcommand{\bbn}{\mathbb N}

\newcommand{\indset}[1]{\textsc{#1}}

\newcommand{\MDS}{\mbox{MDS}}
\newcommand{\TDDP}{\mbox{TDDP}}

\newcommand{\lneg}{^{\neg}}

\newcommand{\rtd}[3] {\mbox{RTD}_{#1}^{#2}(#3)}
\newcommand{\td}[3] {\mbox{TD}^{#1}(#2, #3)}
\newcommand{\order}[1]{\mbox{ord}(#1)}

\usepackage{lineno}

\begin{document}

\title{Classifying the Arithmetical Complexity of Teaching Models}

\titlerunning{Classifying the Arithmetical Complexity of Teaching Models}

\author{Achilles A. Beros$^1$ \and Ziyuan Gao$^2$ \and Sandra Zilles$^2$}
\authorrunning{A.~Beros, Z.~Gao and S.~Zilles}

\institute{Department of Mathematics\\University of Hawaii 
\\\email{beros@math.hawaii.edu} \and
Department of Computer Science\\University of Regina, Regina, SK, Canada S4S 0A2 
\\\email{\{gao257,zilles\}@cs.uregina.ca}}

\maketitle

\begin{abstract}
This paper classifies the complexity of various teaching models by their
position in the arithmetical hierarchy.  In particular, we
determine the arithmetical complexity of the index sets of
the following classes: (1) the class of uniformly r.e.\ families with 
finite teaching dimension, and (2) the class of uniformly r.e.\ 
families with finite positive recursive teaching dimension witnessed by a 
uniformly r.e.\ teaching sequence.  We also derive the
arithmetical complexity of several other decision problems
in teaching, such as the problem of deciding, given an effective 
coding $\{\cL_0,\cL_1,\cL_2,\ldots\}$ of 
all uniformly r.e.\ families, any $e$ such that $\cL_e = \{L^e_0,L^e_1,\ldots,\}$,
any $i$ and $d$, whether or not the teaching dimension
of $L^e_i$ with respect to $\cL_e$ is upper bounded by $d$.     
\end{abstract}

\section{Introduction}

A fundamental problem in computational learning theory is that of
characterising identifiable classes in a given learning model.  
Consider, for example, Gold's \cite{Gold} model of learning from positive data,
in which a learner is fed piecewise with all the positive examples of an
unknown target language -- often coded as a set of natural numbers
-- in an arbitrary order; as the learner processes the data,
it outputs a sequence of hypotheses that must converge syntactically to a 
correct conjecture.  Of particular interest to inductive inference
theorists is the learnability of classes of 
recursively enumerable (r.e.) languages. Angluin \cite{Ang80} demonstrated 
that a uniformly recursive family is learnable in Gold's model if
and only if it satisfies a certain ``tell-tale" condition. 
As a consequence, the family of nonerasing pattern languages over any fixed alphabet
and the family of regular expressions over $\{0,1\}$
that contain no operators other than concatenation and
Kleene plus are both learnable in the limit.  On the other hand,
even a relatively simple class such as 
one consisting of an infinite set $L$ and all the finite subsets of $L$
cannot be learnt in the limit \cite[Theorem I.8]{Gold}.
Analogous characterisations of learnability have since been
discovered for uniformly r.e.\ families as well as for other learning models
such as behaviourally correct learning \cite{dejong96,baliga99}.

Intuitively, the structural properties of learnable families
seem to be related to the ``descriptive complexity" 
of such families.  By fixing a system of describing families 
of sets, one may wish to compare the relative descriptive 
complexities of families identifiable under different criteria. 
One idea, suggested by computability
theory and the fact that many learnability criteria may
be expressed as first-order formulae, is to analyse the quantifier 
complexity of the formula defining the class of learnable families.  
In other words, one may measure the descriptive complexity of identifiable 
classes that are first-order definable by determining the position
of their corresponding index sets in the \emph{arithmetical hierarchy}. 
This approach to measuring the complexity of learnable classes was taken 
by Brandt \cite{brandt86}, Klette \cite{klette76}, and later Beros 
\cite{beros-JSL14}.
More specifically, Brandt \cite[Corollary 1]{brandt86} 
showed that every identifiable class of partial-recursive functions is
contained in another identifiable class with an index
set that is in $\Sigma_3 \cap \Pi_3$, while Beros \cite{beros-JSL14}
established the arithmetical complexity of index sets of uniformly
r.e.\ families learnable under different criteria.  

The purpose of the present work is to determine the arithmetical complexity
of various decision problems in \emph{algorithmic teaching}.
Teaching may be viewed as a natural counterpart of learning,
where the goal is to find a sample efficient learning \emph{and}
teaching protocol that guarantees learning success.  Thus,
in contrast to a learning scenario where the learner has
to guess a target concept based on labelled examples from
a truthful but arbitrary (possibly even adversarial) source,
the learner in a cooperative teaching-learning model
is presented with a sample of labelled examples carefully chosen
by the teacher, and it decodes the sample according to some
pre-agreed protocol.  We say that a family is ``teachable" 
in a model if and only if the associated teaching parameter -- 
such as the \emph{teaching dimension} \cite{GoldmanK95},
the \emph{extended teaching dimension} \cite{hegedus95} or 
the \emph{recursive teaching dimension} \cite{ZillesLHZ11} -- 
of the family is finite.  Due to the ubiquity of
numberable families of r.e.\ sets in theoretical computer
science and the naturalness of such families, 
our work will focus on the class of uniformly r.e.\
families.  Our main results 
classify the arithmetical complexity of index sets of uniformly 
r.e.\ families that are teachable under the teaching dimension
model and a few variants of the recursive teaching dimension model.

From the viewpoint of computability theory, our work provides 
a host of natural examples of complete sets, thus
supporting Rogers' view that many ``arithmetical sets
with intuitively simple definitions $\ldots$ have 
proved to be $\Sigma_n^0$-complete or $\Pi_n^0$-complete (for some $n$)"
\cite[p. 330]{Rog}.  From the viewpoint of computational
learning theory, our results shed light on the recursion-theoretic 
structural properties of the classes of uniformly r.e.\ families that are 
teachable in some well-studied models.  

\vspace*{-.1in}

\section{Preliminaries}

The notation and terminology from computability theory adopted in this
paper follow in general the book of Rogers \cite{Rog}.

$\forall^{\infty}x$ denotes ``for almost
every $x$", $\exists^{\infty}x$ denotes ``for infinitely many
$x$" and $\exists!x$ denotes ``for exactly one $x$". 
$\mathbb{N}$ denotes the set of natural numbers, $\{0,1,2,\ldots\}$,
and $\mathcal{R}$ (= $2^{\natnum}$) denotes the power
set of $\natnum$.  For any function $f$,
$\mbox{ran}(f)$ denotes the range of $f$.
For any set $A$, $\mathbbm{1}_A$ will denote
the characteristic function of $A$, that is, $\mathbbm{1}_A(x) = 1$
if $x \in A$, and $\mathbbm{1}_A(x) = 0$ if $x \notin A$. 
For any sets $A$ and $B$, $A \times B = \{\spn{a,b}: a \in A \wedge
b \in B\}$ and $A\oplus B$, the \emph{join of $A$ and $B$}, 
is the set $\{2x:x\in A\}\cup \{2y+1:y\in B\}$.
Analogously, for any class $\{A_i: i \in \natnum\}$ of sets, 
$\bigoplus_{i\in\natnum}A_i$ is the infinite join of the $A_i$'s.
For any set $A$, $A^*$ denotes the set of all finite sequences
of elements of $A$. 
Given a sequence of families, $\{\mathcal F_i\}_{i\in\bbn}$, 
we define $\bigsqcup_{i\in\bbn} \mathcal F_i = \bigcup_{i\in\bbn}$ 
$\big\{ F\oplus\{i\} : F \in \mathcal F_i \big\}$.  If there are 
only finitely many families, $\mathcal F_0, \ldots, \mathcal F_n$, 
we denote this by $\mathcal F_0 \sqcup \ldots \sqcup \mathcal F_n$.  
We call this the \emph{disjoint union} of
$\mathcal{F}_1,\ldots,\mathcal{F}_n$.

For any $\sigma,\tau\in \{0,1\}^*, \sigma \preceq
\tau$ if and only if $\sigma$ is a prefix of $\tau$,
$\sigma \prec \tau$ if and only if $\sigma$ is a proper prefix of
$\tau$, and $\sigma(n)$ denotes the element in the $n$th position of
$\sigma$, starting from $n=0$.  

Let $W_0,W_1,W_2,\ldots$ be an acceptable numbering of
all r.e.\ sets, and let $D_0,D_1,$ $D_2,\ldots$ be a
canonical numbering of all finite sets such that
$D_0 = \emptyset$ and for any pairwise distinct
numbers $x_1,\ldots,x_n$, $D_{2^{x_1}+\ldots+2^{x_n}}
= \{x_1,\ldots,x_n\}$.  For all $e$ and $j$, define
$L^e_j = \{x: \spn{j,x} \in W_e\}$ and
$\mathcal{L}_e = \{ L^e_j: j \in \natnum\}$.
$L^e_j$ is the $j$th column of $W_e$. 
$L^e_{j,s}$ denotes the $s$th approximation of $L^e_j$ which, 
without loss of generality, we assume is a subset of  
$\{0,\ldots,s\}$.  Note that $\{ \mathcal{L}_e: e\in\natnum\}$
is the class of all uniformly r.e.\ (u.r.e.) families,
each of which is encoded as an r.e.\ set.
Let \indset{coinf} denote the index
set of the class of all coinfinite r.e.\ sets, and let
\indset{cof} denote the index set of the class of all cofinite
r.e.\ sets.  Let \indset{inf} denote the index set of the class of all infinite 
r.e.\ sets and \indset{fin} denote the index set of the class of all finite sets.    

\begin{defn}\cite{Rog}
A set $A \subseteq \natnum$ is in $\Sigma_0 (=\Pi_0 = \Delta_0)$ iff $A$ is recursive.
$A$ is in $\Sigma^0_n$ iff there is a recursive relation $R$ such that
\begin{equation}\label{defn:sigman}
x \in B \leftrightarrow (\exists y_1)(\forall y_2)\ldots(Q_n y_n)
R(x,y_1,y_2,\ldots,y_n)
\end{equation}
where $Q_n = \forall$ if $n$ is even and $Q_n = \exists$ if $n$ is odd.
A set $A \subseteq \natnum$ is in $\Pi^0_n$ iff its complement
$\overline{A}$ is in $\Sigma^0_n$. 
$(\exists y_1)(\forall y_2) \ldots (Q_n y_n)$ is known as a $\Sigma^0_n$ prefix;
$(\forall y_1)(\exists y_2) \ldots (Q_n,y_n)$, where $Q_n = \exists$
if $n$ is even and $Q_n = \forall$ if $n$ is odd, is known as a $\Pi_n^0$ prefix.
The formula on the right-hand side of (\ref{defn:sigman})
is called a \emph{$\Sigma_n^0$ formula} and its negation 
is called a \emph{$\Pi_n^0$ formula}. 
A set $A$ is in $\Delta_n^0$ iff $A$ is in $\Sigma_n^0$ and
$A$ is in $\Pi_n^0$.  Sets in $\Sigma_n^0$ ($\Pi_n^0$, $\Delta_n^0$)
are known as \emph{$\Sigma_n^0$ sets} (\emph{$\Pi_n^0$ sets, $\Delta_n^0$ sets}).
For any $n \geq 1$, a set $A$ is \emph{$\Sigma_n^0$-hard} 
(\emph{$\Pi_n^0$-hard}) iff
every $\Sigma_n^0$ ($\Pi_n^0$) set $B$ is many-one reducible to it, that is,
there exists a recursive function $f$ such that
$x \in B \leftrightarrow f(x) \in A$. 
$A$ is \emph{$\Sigma_n^0$-complete} (\emph{$\Pi_n^0$-complete}) iff 
$A$ is definable with a $\Sigma_n^0$ ($\Pi_n^0$) formula and
$A$ is $\Sigma_n^0$-hard (\emph{$\Pi_n^0$-hard}).
\end{defn}

\noindent The following proposition collects several useful equivalent 
forms of $\Sigma_n^0$ or $\Pi_n^0$ formulas (for any $n$).

\begin{prop}\label{prop:sigmaequivform}
\begin{enumerate}[{\sc (i)}]
\item For every $\Sigma_{n+1}^0$ set $A$, there is a
$\Pi_n^0$ predicate $P$ such that for all $x$,
$$
x \in A \leftrightarrow (\forall^{\infty}a)P(a,x) \leftrightarrow (\exists a)P(a,x).
$$ 
\item For every $\Sigma_{n+1}^0$ set $B$, there is a
$\Pi_n^0$ predicate $Q$ such that for all $x$,
$$
x \in B \leftrightarrow (\exists! a)Q(a,x) \leftrightarrow (\exists a)Q(a,x).
$$ 
\end{enumerate}
\end{prop}

\section{Teaching}

Goldman and Kearns \cite{GoldmanK95} introduced a variant
of the on-line learning model in which a helpful teacher selects
the instances presented to the learner.  They considered
a combinatorial measure of complexity called the \emph{teaching dimension},
which is the mininum number of labelled examples required
for any consistent learner to uniquely identify any target concept
from the class.  

Let $\mathcal{L}$ be a family of subsets of
$\mathbb{N}$.  Let $L \in \mathcal{L}$
and $T$ be a subset of $\mathbb{N} \times \{+,-\}$.
Furthermore, let $T^+ = \{ n : (n,+) \in 
T \}$, $T^- = \{ n : (n,-) \in T \}$ and $X(T) = T^+ \cup 
T^-$. 
A subset $L$ of $\natnum$ is
said to be \emph{consistent} with $T$ iff
$T^+ \subseteq L$ and $T^- \cap L = \emptyset$.
$T$ is a \emph{teaching set} for $L$ with respect to $\mathcal{L}$
iff $T$ is consistent with $L$ and for all $L'\in\mathcal{L}\sm\{L\}$,
$T$ is not consistent with $L'$. 
If $T$ is a teaching set for $L$ with respect to $\mathcal{L}$,
then $X(T)$ is known as a \emph{distinguishing set} for $L$
with respect to $\mathcal{L}$.  
Every element of $\mathbb{N} \times \{+,-\}$
is known as a \emph{labelled example}.

\begin{defn}\cite{GoldmanK95,ShinoharaM91}
Let $\mathcal{L}$ be any family of subsets  
of $\mathbb{N}$. Let $L\in\mathcal{L}$ be given. 
The size of a smallest teaching set for $L$ with respect to $\mathcal{L}$
is called the \emph{teaching dimension of $L$ with respect to 
$\mathcal{L}$}, denoted by $\TD(L,\mathcal{L})$. 
The \emph{teaching dimension of $\mathcal{L}$} is 
defined as $\sup\{\TD(L,\mathcal{L}): 
L\in\mathcal{L}\}$ and is denoted by $\TD(\mathcal{L})$.
If there is a teaching set for $L$ with respect to $\mathcal{L}$
that consists of only positive examples, 
then the \emph{positive teaching dimension of $L$ with respect to
$\mathcal{L}$} is defined to be the smallest possible size
of such a set, and is denoted by $\TD^+(L,\mathcal{L})$.
If there is no teaching set for $L$ w.r.t. $\mathcal{L}$
that consists of only positive examples,
then $\TD^+(L,\mathcal{L})$ is defined to be $\infty$.
A teaching set for $L$ with respect to $\cL$ that consists
of only positive examples is known as a \emph{positive
teaching set} for $L$ with respect to $\cL$.
The \emph{positive teaching dimension of $\mathcal{L}$} is
defined as $\sup\{\TD^+(L,\mathcal{L}): L \in \mathcal{L}\}$.
\end{defn}

\noindent Another complexity parameter recently studied in computational 
learning theory is the recursive teaching dimension. It refers 
to the maximum size of teaching sets in a series of nested subfamilies 
of the family.

\begin{defn} (Based on \cite{ZillesLHZ11,zeinab})\label{rtddef}
Let $\mathcal{L}$ be any family of 
subsets of $\mathbb{N}$.
A \emph{teaching sequence for $\mathcal{L}$} 
is any sequence $\TS = ((\mathcal{F}_0,d_0), (\mathcal{F}_1,d_1),\linebreak[3]\ldots)$ 
where (i)~the families $\mathcal{F}_i$ form a partition 
of $\mathcal{L}$ with each $\mathcal{F}_i$ nonempty, 
and (ii)~$d_i = \sup\{\TD(L,\mathcal{L}\setminus\bigcup_{0\le j<i}\mathcal{F}_j):
L \in \mathcal{F}_i\}$ for all $i$.  $\sup\{d_i: i
\in\mathbb{N}\}$ is called the \emph{order of $\TS$}, and is denoted by $ord(\TS)$.  
The \emph{recursive teaching dimension of $\mathcal{L}$} is defined
as $\inf\{ord(\TS):\mbox{$\TS$ is a}$ $\mbox{teaching sequence for $\mathcal{L}$}\}$ 
and is denoted by $\RTD(\mathcal{L})$.
\end{defn}

\noindent We shall also briefly consider the \emph{extended teaching dimension} (\XTD)
of a class.  This parameter may be viewed as a generalisation of the
teaching dimension; it expresses the complexity of unique
specification with respect to a concept class $\cC$ for
\emph{every} concept (not just members of $\cC$) over a given instance space $X$.
As Heged\"{u}s \cite{hegedus95} showed, the extended teaching dimension
of a concept class $\cC$ is closely related to the query complexity of learning $\cC$.

\begin{defn}\cite{hegedus95}
Let $\cL$ be a family of subsets of $\natnum$,
and let $L$ be a subset of $\natnum$.
A set $S \subseteq \natnum$ is a \emph{specifying set for $L$
with respect to $\cL$} iff there is at most one concept $L'$
in $\cL$ such that $L \cap S = L' \cap S$. 
Define the \emph{extended teaching dimension} (\XTD) of $\cL$
as $\inf\{d:\mbox{for every set $L \subseteq \natnum$ there exists
an at most}$ $\mbox{$d$-element specifying set with respect to $\cL$}\}$.

A set $S \subseteq L$ is a \emph{positive specifying set for $L$
with respect to $\cL$} iff there is at most one concept
in $\cL$ that contains $S$.
Define the \emph{positive extended teaching dimension} ($\XTD^+$) of $\cL$
as $\inf\{d:\mbox{for every set $\emptyset \neq L \subseteq \natnum$ there exists
an at most}$ $\mbox{$d$-element positive specifying set with respect to $\cL$}\}$.
If there is a nonempty set $L$ that does not have a positive specifying
set w.r.t\ $\cL$, define $\XTD^+(\cL) = \infty$.   
\end{defn}

\noindent The next series of definitions will introduce various subsets 
of $\natnum$, each of which is a set of codes for u.r.e.\ families that
satisfy some notion of teachability.

\begin{enumerate}[{\sc (i)}]
\item $I_{TD}^{\forall} = \{ e: (\forall L\in\mathcal{L}_e)
[\TD(L,\mathcal{L}_e) < \infty] \}$.
\item $I_{TD} = \{ e: \TD(\mathcal{L}_e) < \infty\}$.
\item $I_{TD}^{\forall^{\infty}} = \{ e: (\forall^{\infty} L \in \mathcal{L}_e)
[\TD(L,\mathcal{L}_e) < \infty]\}$.
\item $I_{TD^+}^{\forall} = \{e: (\forall L \in \mathcal{L}_j)[\TD^+(L,\mathcal{L}_e) < \infty]\}$.
\item $I_{TD^+} = \{ e: \TD^+(\mathcal{L}_e) < \infty\}$.
\item $I_{TD^+}^{\forall^{\infty}} = \{ e: (\forall^{\infty} L \in \mathcal{L}_e)
[\TD^+(L,\mathcal{L}_e) < \infty]\}$.
\item $I_{XTD} = \{ e: \XTD(\mathcal{L}_e) < \infty\}$.
\item $I_{XTD^+} = \{ e: \XTD^+(\mathcal{L}_e) < \infty\}$.
\end{enumerate}

\noindent Owing to space constraints, many proofs will be omitted or
sketched.

\section{Teaching Dimension}

In this section we study the arithmetical
complexity of the class of u.r.e.\ families 
with finite teaching dimension;
several related decision problems 
will also be considered. 

Before proceeding with the main theorems on the arithmetical
complexity of the teaching dimension model and its variants, 
a series of preparatory results will be presented.
Theorem \ref{thm:acds} addresses the question: 
how hard (arithmetically) is it to determine whether 
or not, given $e\in\natnum$ and a finite set $D$, 
$D$ can distinguish an r.e.\ set 
$L_j^e \in \cL_e$ from the other members of $\cL_e$?

\begin{defn}
$\DS := \{ \langle e,x,u\rangle: (\forall y)
[L^e_x \neq L^e_y \rightarrow L^e_x \cap D_u \neq L^e_y \cap D_u]\}$.\footnote{\DS\ stands for
``distinguishing set."}
\end{defn} 

\begin{thm}\label{thm:acds}
\DS\ is $\Pi_2^0$-complete.
\end{thm}

\proof
By the definition of \DS, $\langle e,x,u\rangle \in \DS \leftrightarrow 
(\forall y)(\forall t)(\exists v > t)[[L^e_{x,v} \cap D_u \neq L^e_{y,v} \cap D_u]
\vee (\forall p)(\forall a)(\exists b > a)[\charf_{L^e_{x,b}}(p) = \charf_{L^e_{y,b}}(p)]]$.
Thus \DS\ has a $\Pi_2^0$ description.
Now, since \indset{inf} is $\Pi_2^0$-complete \cite{Rog},
it suffices to show that \indset{inf} is many-one reducible to
\DS.  Let $g$ be a one-one recursive function with
$$
W_{g(e,i)} =  \begin{cases}
\natnum & \mbox{if $i = 0 \vee (\exists j > i)[j \in W_e]$;} \\
\{0,1,\ldots,i\} & \mbox{otherwise.}\end{cases}  
$$
Let $f$ be a recursive function such that 
$\cL_{f(e)} = \{W_{g(e,i)}: i \in \natnum\}$
and $L^{f(e)}_0 = W_{g(e,0)}$.
Recall that $D_0 = \emptyset$.
It is readily verified that $e \in \indset{inf} \leftrightarrow
\langle f(e),0,0\rangle \in \DS$.~\qed

\medskip
\noindent
The expectation that the arithmetical
complexity of determining if a finite $D$
is a \emph{smallest} possible distinguishing
set for some $W_x$ belonging to $\mathcal{L}_e$ is
at most one level above that of \DS\ is confirmed  
by Theorem~\ref{mdscomplete}.

\begin{defn}
$\MDS \!:=\! \{ \langle e,x,u \rangle \in \DS\! : \!(\forall u')[|D_{u'}| < |D_u|
\rightarrow \langle e,x,u' \rangle \notin \DS]\}$.\footnote{\MDS\ stands for
``minimal distinguishing set."}
\end{defn}

\begin{thm}\label{mdscomplete}
\MDS\ is $\Pi_3^0$-complete.
\end{thm}

\proof (Sketch.)
By the definition of \MDS,
\[ \langle e,x,u \rangle \in \MDS \leftrightarrow
\langle e,x,u \rangle \in \DS \wedge (\forall u')[(|D_{u'}| \geq |D_u|)
\vee \langle e,x,u' \rangle \notin \DS]. \]
By Theorem \ref{thm:acds}, \DS\ has a $\Pi_2^0$ description and
$\overline{\DS}$ has a $\Sigma_2^0$ description.
Thus \MDS\ has a $\Pi_3^0$ description.
We omit the proof that \MDS\ is $\Pi_3^0$-hard.~\qed

\medskip
\noindent
Another problem of interest is the complexity of determining
whether or not the teaching dimension of some $W_x$ w.r.t. a class 
$\mathcal{L}_e$ is upper-bounded by a given number $d$.
For $d = 0$, this problem is just as
hard as \DS\ (see Proposition \ref{prop:actddp0}); for $d > 0$, however, the complexity of the problem 
is exactly one level above that of \DS\ (see Theorem \ref{thm:actddp}).  We omit the proofs.

\begin{defn}
$\TDDP := \{ \langle e,x,d \rangle: d \geq 1 \wedge (\exists u)[|D_u| \leq d \wedge
\langle e,x,u \rangle \in \DS]\}$.\footnote{\TDDP\ stands for ``Teaching dimension
decision problem."}
\end{defn}  

\begin{prop}\label{prop:actddp0}
$\{\langle e,x\rangle: \mathcal{L}_e = \{L^e_x\}\}$ is $\Pi_2^0$-complete.
\end{prop}

\begin{thm}\label{thm:actddp}
\TDDP\ is $\Sigma_3^0$-complete.
\end{thm}

\medskip
\noindent
Our first main result states that the class of all  
u.r.e.\ families $\cL$ such that any finite subclass $\cL' \subseteq \cL$ has
finite teaching dimension with respect to $\cL$ is $\Pi_4$-complete. 
 
\begin{thm}\label{t1complete}
$I_{TD}^{\forall}$ is $\Pi_4^0$-complete.
\end{thm}

\proof
First, note 
the following equivalent
conditions:
\[ e \in I_{TD}^{\forall} \leftrightarrow (\forall i)[\TD(L_i^e,\cL_e) < \infty)] \leftrightarrow
(\forall i)(\exists u)[\spn{e,i,u} \in \DS].  \]
By Theorem \ref{thm:acds}, $\spn{e,i,u} \in \DS$ may be
expressed as a $\Pi^0_2$ predicate, so
$I_{TD}^{\forall}$ is $\Pi_4^0$.

Now consider any $\Pi_4^0$ unary predicate $P(e)$;
$P(e)$ is of the form $(\forall x)[Q(e,x)]$, where $Q$
is a $\Sigma_3^0$ predicate.  Since \indset{cof} is $\Sigma_3^0$-complete \cite{Rog},
there is a recursive function $g(e,x)$ such that
$P(e) \leftrightarrow (\forall x)[Q(e,x)] \leftrightarrow
(\forall x)[g(e,x) \in \indset{cof}]$ holds.  
For each triple $\langle e,x,i\rangle$, define 
$$
H_{\langle e,x,i\rangle} = \begin{cases}
\{ \langle e,x\rangle \} \oplus (W_{g(e,x)} \cup \{i\}) & \mbox{if $i > 0$;} \\
\{ \langle e,x\rangle \} \oplus W_{g(e,x)} & \mbox{if $i = 0$}\end{cases}
$$   
Let $h$ be a recursive function such that for all $e$, 
$\cL_{h(e)} = \{H_{\spn{e,x,i}}: x,i\in\natnum\}$.
\begin{description} 
\item[Case (i):] $P(e)$ holds.  Then for
all $x$, $W_{g(e,x)}$ is cofinite.  Thus
for all $x$ and each $i > 0$ such that $i \notin W_{g(e,x)}$, 
$H_{\langle e,x,i\rangle}$ has
the teaching set $\{(2\langle e,x\rangle,+),(2i+1,+)\}$
with respect to $\mathcal{L}_{h(e)}$.  
Furthermore, for all $x$ and each
$i$ such that either $i \neq 0 \wedge i \in W_{g(e,x)}$ or $i = 0$,  
$H_{\langle e,x,i\rangle}$ has the teaching set
$\{ (2\langle e,x\rangle,+) \} \cup \{(2j+1,-): j\notin W_{g(e,x)}\wedge
j > 0\}$ with respect to $\mathcal{L}_{h(e)}$.
Therefore $\TD(H_{\spn{e,x,i}},$ $\mathcal{L}_{h(e)}) < \infty$
for every pair $\langle x,i\rangle$, so that $h(e)\in I_{TD}^{\forall}$.

\item[Case (ii):] $\lneg P(e)$ holds.  Then $W_{g(e,x)}$ is coinfinite
for some $x$.  Fix such an $x$.  Then $\mathcal{L}_{h(e)}$ contains
$\mathcal{L}' = \{H_{\langle e,x,i\rangle}: i\in\natnum\}$.
Furthermore, for each positive $i\notin W_{g(e,x)}$,
since $\{\langle e,x\rangle \} \oplus (W_{g(e,x)} \cup \{i\}) \in \mathcal{L}'$,
any teaching set for $H_{\langle e,x,0\rangle}$
w.r.t.\ $\mathcal{L}_{h(e)}$ must contain $(2i+1,-)$.  Hence 
$\TD(H_{\spn{e,x,0}},\mathcal{L}_{h(e)}) = \infty$,
so that $h(e) \notin I_{TD}^{\forall}$.
\end{description}
Thus $I_{TD}^{\forall}$ is $\Pi_4^0$-complete.~\qed

\medskip
\noindent Extending $I_{TD}^{\forall}$ to include u.r.e.\ families
$\cL$ for which there is a \emph{cofinite} subclass $\cL' \subseteq \cL$
belonging to $I_{TD}^{\forall}$ increases the arithmetical complexity of $I_{TD}^{\forall}$ to $\Sigma_5^0$.

\begin{thm}\label{t3complete}
$I_{TD}^{\forall^{\infty}}$ is $\Sigma_5^0$-complete.
\end{thm}

\proof (Sketch.)
For any $e,s$, let $g$ be a recursive function such that
$\cL_{g(e,s)} = \{L^e_i: i > s\}$.  Note that the 
expression for $I_{TD}^{\forall^{\infty}}$ can be re-written
as $e \in I_{TD}^{\forall^{\infty}} \leftrightarrow (\exists t)(\forall s > t)[g(e,s) \in I_{TD}^{\forall}].$
Since $\mathcal{L}_e \sm \cL_{g(e,t)}$ is finite,
it follows that $e \in I_{TD}^{\forall^{\infty}}$.  By Theorem \ref{t1complete},
the predicate $g(e,s) \in I_{TD}^{\forall}$ has a $\Pi_4^0$
description, so that $I_{TD}^{\forall^{\infty}}$ is definable with a $\Sigma_5^0$
predicate.

Now let $P$ be any $\Sigma_5^0$ predicate.  By Proposition
\ref{prop:sigmaequivform} and the $\Sigma_3^0$-completeness of 
\indset{cof} \cite{Rog}, there is a recursive function $h$ such
that $P(e) \leftrightarrow (\forall^{\infty}a)(\forall b)[h(e,a,b) \in \indset{cof}]
\leftrightarrow (\exists a)(\forall b)[h(e,a,b) \in \indset{cof}]$  
and
$(\exists b)[h(e,a,b) \in \indset{coinf}] \leftrightarrow (\exists ! b)
[h(e,a,b)$ $\in \indset{coinf}].$
Now let $g$ be a one-one recursive function such that
$$
W_{g(e,a,b,i)} = \begin{cases}
\{ \langle e,a,b\rangle \} \oplus (W_{h(e,a,b)} \cup \{i\}) & \mbox{if $i > 0$;} \\
\{ \langle e,a,b\rangle \} \oplus W_{h(e,a,b)} & \mbox{if $i = 0$.}\end{cases}
$$ 
Let $f$ be a one-one recursive function such that
$\cL_{f(e)} = \{ W_{g(e,a,b,i)}:a,b,i\in\natnum\}$.
Note that for all $a,b,i\in\natnum$, $\TD(W_{g(e,a,b,i)},\mathcal{L}_{f(e)}) < \infty
\leftrightarrow h(e,a,b) \in \indset{cof}$. 
One can show as in the proof of Theorem \ref{t1complete}
that $\TD(L,\cL_{f(e)}) < \infty$ holds for almost all $L \in \cL_{f(e)}$
iff $P(e)$ is true.~\qed

\medskip
\noindent The next theorem shows that the index set of 
the class consisting of all u.r.e.\ families
with finite teaching dimension is $\Sigma_5^0$-complete. 

\begin{thm}\label{thm:td}
$I_{TD}$ is $\Sigma_5^0$-complete.
\end{thm} 

\proof (Sketch.)
From $\TD(\cL_e) < \infty \leftrightarrow (\exists a)(\forall b)[\spn{e,b,a} \in \TDDP]$
and the fact that $\TDDP$ is $\Sigma_3^0$-complete by Theorem \ref{thm:actddp}
we have that $I_{TD}$ is $\Sigma_5^0$.

To prove that $I_{TD}$ is $\Sigma_5^0$-hard, consider any $\Sigma_5^0$ predicate
$P(e)$.  There is a binary recursive function $g$ such that
$P(e) \rightarrow (\exists a)(\forall b)[g(e,a,b) \in \indset{cof}]$  
and $\lneg P(e) \rightarrow (\forall a)(\forall^{\infty}b)[g(e,a,b) \in \indset{coinf}].$
Now fix $e,b \in \natnum$.  For each $a$, let $\{H_{\spn{a,0}},H_{\spn{a,1}},H_{\spn{a,2}},\ldots\}$ be a
numbering of the union of two u.r.e.\ families
$\{C_{\spn{a,0}},C_{\spn{a,1}},C_{\spn{a,2}},\ldots\}$ and $\{L_{\spn{a,0}},L_{\spn{a,1}},L_{\spn{a,2}},\ldots\}$, 
which are defined as follows.  (For simplicity, the notation for dependence
on $e$ and $b$ is dropped.)

\begin{enumerate}
\item $\{C_{\spn{a,0}},C_{\spn{a,1}},C_{\spn{a,2}},\ldots\}$ is
a numbering of $\{X \subseteq \natnum: |\natnum \sm X| = a\}$.
\item Let $E_{\spn{a,0}},E_{\spn{a,1}},E_{\spn{a,2}},\ldots$ be
a one-one enumeration of $\{X \subseteq \natnum: |\natnum \sm X| < \infty \wedge |\natnum \sm X| \neq a\}$.
Let $f$ be a recursive function such that for all $n,s\in\natnum$,
$f(n,s)$ is the $(n+1)$st element of $\natnum \sm W_{g(e,a,b),s}$. 
For all $n,s \in \natnum$, define
$$
L_{\spn{a,\spn{n,s}}} = \begin{cases}
E_{\spn{a,n}} & \mbox{if $(\forall t \geq s)[f(n,s) = f(n,t)]$;} \\
\natnum & \mbox{if $(\exists t > s)[f(n,s) \neq f(n,t)]$.}\end{cases}
$$
\end{enumerate}

Note that $\{H_{\spn{a,0}},H_{\spn{a,1}},H_{\spn{a,2}},\ldots\} \neq \emptyset$.
Now construct a u.r.e.\ family $\{G_{\spn{e,b,0}},G_{\spn{e,b,1}},\ldots\}$
with the following properties:

\begin{enumerate}[{\sc (i)}]
\item For all $s$, $G_{\spn{e,b,s}}$ is of the form $\{b\} \oplus \{s\}
\oplus \displaystyle\bigoplus\nolimits_{j \in\natnum} H_{\spn{j,i_j}}$.
\item For every nonempty finite set $\{H_{\spn{c_0,d_0}},\ldots,H_{\spn{c_k,d_k}}\}$
with $c_0 < \ldots < c_k$, there is at least one $t$ for which 
$G_{\spn{e,b,t}} = \{b\} \oplus \{t\} \oplus 
\displaystyle\bigoplus\nolimits_{i\in\natnum}A_i$, 
where $A_{c_i} = H_{\spn{c_i,d_i}}$ for all $i \in \{0,\ldots,k\}$
and $A_i \in \{H_{\spn{i,j}}:j\in\natnum\}$ for all $i \notin \{c_0,\ldots,c_k\}$.
\item For every $t$ such that $G_{\spn{e,b,t}} = \{b\} \oplus \{t\} \oplus
\displaystyle\bigoplus\nolimits_{j \in\natnum} H_{\spn{j,i_j}}$, 
there are infinitely many $t' \neq t$ such that
$G_{\spn{e,b,t'}} = \{b\} \oplus \{t'\} 
\oplus \displaystyle\bigoplus\nolimits_{j \in\natnum} H_{\spn{j,i_j}}$.
\end{enumerate}

The family $\{G_{\spn{e,i,j}}: i,j \in \natnum\}$ may be thought
of as an infinite matrix $M$ in which each row represents a 
set parametrised by $g(e,a,b)$ for a fixed $b$ and
$a$ ranging over $\natnum$.  
Furthermore, if there exists an $a$ such that 
$W_{g(e,a,b)}$ is cofinite for all $b$, then
the $a$th column of $M$ contains all cofinite 
sets with complement of size $a$ plus a finite
number of other cofinite sets; if no such $a$ exists
then every column of $M$ contains all cofinite sets.
Let $h$ be a recursive function with $\cL_{h(e)} = 
\{\{b\} \oplus \displaystyle\bigoplus\nolimits_{i \in\natnum} \emptyset: b \in \natnum\}
\cup \{G_{\spn{e,b,s}}: b,s \in \natnum\}$. Showing that
$P(e)$  iff $h(e) \in I_{TD}$ proves $I_{TD}$ to be $\Sigma_5^0$-complete.  
~\qed

\medskip
\noindent To conclude our discussion on the general teaching dimension, we 
demonstrate 
that the criterion for a u.r.e.\ family to have finite 
extended teaching dimension 
is so stringent that only 
finite families have a finite \XTD.

\begin{thm}\label{thm:xtd}
$I_{XTD}$ is $\Sigma_3$-complete.
\end{thm}   

\proof
We show $\XTD(\cL_e) < \infty \leftrightarrow
|\cL_e| < \infty$.  First, suppose $\cL_e = \{L_1,\ldots,L_k\}$.
As $\cL_e$ is finite, $\TD(L_i,\cL_e) \leq k-1$ for all 
$i \in \{1,\ldots,k\}$.  Consider any $L \notin \cL_e$.
For each $i \in \{1,\ldots,k\}$, fix $y_i$ with $\mathbbm{1}_L(y_i) \neq \mathbbm{1}_{L_i}(y_i)$.
Then $\{ (y_i,+): 1 \leq i \leq k \wedge y_i \in L\} \cup 
\{(y_j,-): 1 \leq j \leq k \wedge y_j \notin L\}$
is a specifying set for $L$ with respect to $\cL_e$ of size $k$.

Second, suppose $|\cL_e| = \infty$.  Let $T = \{\sigma\in\{0,1\}^*:
(\exists^{\infty}L \in \cL_e)(\forall x < |\sigma|)[\sigma(x) = \mathbbm{1}_L(x)]\}$.      
Note that $|\cL_e| = \infty$ implies $T$ is an infinite binary tree.  
Thus by K\"{o}nig's lemma, $T$ has at least one infinite branch, say $B$.
Then for all $n \in \natnum$, there exist infinitely many
$L \in \cL_e$ such that $\mathbbm{1}_B(x) = \mathbbm{1}_L(x)$
for all $x < n$.  Therefore $B$ has no finite specifying set
with respect to $\cL_e$ and so $\XTD(\cL_e) = \infty$.
Consequently, $\XTD(\cL_e) < \infty \leftrightarrow
|\cL_e| < \infty \leftrightarrow (\exists a)(\forall b)(\exists c \leq a)
(\forall x)(\forall s)(\exists t > s)[\charf_{L_{b,t}^e}(x) = \charf_{L^e_{c,t}}(x)]$;
as any two quantifiers, at least one of which is bounded,
may be permuted, it follows that the last expression
is equivalent to a $\Sigma_3^0$ formula.
To show that $I_{XTD}$ is $\Sigma_3^0$-hard, consider any
$\Sigma_3^0$ predicate $P$, and let $g$ be
a recursive function such that
\[ P(e) \leftrightarrow g(e) \in \indset{cof}.  \]   
Let $f$ be a recursive function with
$\cL_{f(e)}$ equal to $\{W_{g(e)} \cup D: \natnum \supseteq |D| < \infty\}$,
the class of all sets consisting of the union of $W_{g(e)}$
and a finite set of natural numbers.
Then
\[ P(e) \leftrightarrow g(e) \in \indset{cof} \leftrightarrow
|\cL_{f(e)}| < \infty \leftrightarrow \XTD(\cL_{f(e)}) < \infty,\]
and so $I_{XTD}$ is indeed $\Sigma_3^0$-complete.~\qed

\section{Positive Teaching Dimension}
 
We now consider the analogues of the results in the preceding
section for the positive teaching dimension model.  In studying
the process of child language acquisition, Pinker \cite[p. 226]{pinker79}
points to evidence in prior research that children are seldom ``corrected
when they speak ungrammatically", and ``when they are corrected they
take little notice".  It seems likely, therefore, that children learn
languages mainly from positive examples.  Thus, as a model for
child language acquisition, the positive teaching dimension model 
is probably closer to reality than the general teaching dimension model 
in which negative examples are allowed.       
The next two results are the analogues of Theorems \ref{t1complete} 
and \ref{thm:td} for the positive teaching dimension model.
It is noteworthy that $I_{TD}^{\forall}$ and $I_{TD^+}^{\forall}$ have equal arithmetical complexity; 
that is to say, restricting the teaching sets of each $L \in \cL_e$ with 
$e \in I_{TD}^{\forall}$ to positive teaching sets has no overall effect on the arithmetical complexity of $I_{TD}^{\forall}$.  

\begin{thm}\label{tdpluseachfinite}
$I_{TD^+}^{\forall}$ is $\Pi_4^0$-complete.
\end{thm}

\proof (Sketch.)
Observe that
\[ (\forall L \in \mathcal{L}_e)[\TD^+(L,\mathcal{L}_e) < \infty] \leftrightarrow
(\forall i)(\exists u)[(\exists s)[D_u \subseteq L_{i,s}^e] \] 
\[ \wedge
(\forall j)[(\forall x)(\forall s)(\exists t>s)[\charf_{L^e_{j,t}}(x) = \charf_{L^e_{i,t}}(x)]
\vee \forall s'[D_u \not\subseteq L^e_{j,s'}]]], \]
Since the right-hand side simplifies to a $\Pi_4^0$ predicate,
we know that $I_{TD^+}^{\forall}$ is $\Pi_4^0$. 

For the proof that $I_{TD^+}^{\forall}$ is $\Pi_4^0$-hard, take any $\Pi_4^0$ predicate
$P$, and let $g$ be a recursive function such that
$P(e) \leftrightarrow (\forall a)[g(e,a) \in \indset{cof}].$
Define a u.r.e.\ family $\cL = \{L_i\}_{i\in\natnum}$ as follows.
(For notational simplicity, the notation for dependence on $e$
is dropped.)  For all $a,i\in\natnum$,
$$
L_{\spn{a,0}} = \{a\} \oplus W_{g(e,a)},
$$
$$
L_{\spn{a,i+1}} = \begin{cases}
\{a\} \oplus W_{g(e,a)} & \mbox{if $i \in W_{g(e,a)}$;} \\
\{a\} \oplus (\{i\} \cup \{x: x < i \wedge x \in W_{g(e,a)}\}) & \mbox{if $i \notin W_{g(e,a)}$.}\end{cases}
$$ 
Let $f$ be a recursive function for which $\cL_{f(e)} = \cL$.
One can show that \linebreak[4]$\TD^+(L,\cL_{f(e)}) < \infty$ holds
for all $L \in \cL_{f(e)}$ iff $P(e)$ is true.~\qed

\begin{thm}\label{tdplusfinite}
$I_{TD^+}$ is $\Sigma_5^0$-complete.
\end{thm}

\proof (Sketch.)
For any $e$, one has 
$\TD^+(\cL_e) < \infty \leftrightarrow (\exists a)(\forall b)[\TD^+(L^e_b,\cL_e) < a]$ and
\[ \TD^+(L^e_b,\cL_e) < a \leftrightarrow (\exists u)(\forall c)[(\exists s')[D_u \subseteq L^e_{b,s'}] 
\wedge |D_u| < a \]
\[ \wedge [(\forall x)(\forall s)(\exists t > s) [\charf_{L^e_b}(x) = \charf_{L^e_c}(x)] \vee (\forall t')[D_u 
\not\subseteq L^e_{c,t'}]]].\] 
Simplifying the last equivalence yields
a $\Sigma_5^0$-predicate for $\TD^+(\cL_e) < \infty$.

The proof that $I_{TD^+}$ is $\Sigma_5^0$-hard is similar to the earlier proof
that $I_{TD}$ is $\Sigma_5^0$-hard (but requires some additional ideas).  
Given any $\Sigma_5^0$ formula $P$, let 
$R$ be a recursive predicate such that
$P(e) \rightarrow (\exists a)(\forall b)(\forall^{\infty}c)(\forall d)(\exists l)
[R(e,a,b,c,d,l)]$ and
$\lneg P(e) \rightarrow (\forall a)(\forall^{\infty}b)(\forall c)(\exists d)(\forall l)$
$[\lneg R(e,a,b,c,d,l)].$
Now fix any $e,b \in \natnum$.  Let $( \langle a^0_0,b^0_0,c^0_0\rangle,\ldots,\langle a^0_{k_0},b^0_{k_0},$ $c^0_{k_0}\rangle)
,$ $(\langle a^1_0,b^1_0,c^1_0 \rangle,\ldots,\langle a^1_{k_1},$ $b^1_{k_1},c^1_{k_1}
\rangle),\ldots$ be a one-one enumeration of all non-empty finite sequences of triples such that
$a^j_0 < \ldots < a^j_{k_j}$ for all $j \in \natnum$. 
Define the set (dropping the notation for dependence on $e$) 
\[ S_b := \{i: (\exists j\in\{0,\ldots,k_i\})[|D_{b^i_j}| \neq a^i_j+1 \wedge
(\exists l)[R(e,a^i_j,b,b^i_j,c^i_j,l)]]\}, \]
using our fixed numbering $D_0,D_1,D_2,\ldots$ of all finite sets.
For each $b$ (with $e$ fixed), construct a u.r.e.\ family $\{G_{\spn{b,-1}}\} \cup 
\bigcup_{s\in\natnum}$ $\{G_{\spn{b,s}}\}$ as follows. 
\[ G_{\spn{b,s}} = \{b\} \oplus S_b \oplus \displaystyle\bigoplus\nolimits_{i \in\natnum}\natnum\
\mbox{if $s = -1$ or $s \in S_b$.} \]
If $s \notin S_b$, set $G_{\spn{b,s}} = \{b\} \oplus (S_b \cup \{s\})
\oplus \displaystyle\bigoplus\nolimits_{i \in\natnum} H_i$, where
$$
H_i =  \begin{cases}
\emptyset & \mbox{if $i \notin \{a^s_0,\ldots,a^s_{k_s}\}$;} \\
D_{b^s_i} & \mbox{if $i \in \{a^s_0,\ldots,a^s_{k_s}\}$.} \end{cases}
$$ 
Let $f$ be a recursive function such that
$\cL_{f(e)} = \bigcup_{b \in \natnum} (\{G_{\spn{b,-1}}\} \cup 
\bigcup_{s\in\natnum}$ $\{G_{\spn{b,s}}\})$ (note again
that the notation for dependence on $e$ in the definition
of $G_{\spn{b,s}}$ has been dropped).
We omit the proof that $P(e)$ holds iff $\TD^+(\cL_{f(e)}) < \infty$. 
~\qed

\begin{thm}\label{tdplusaefinite}
$I_{TD^+}^{\forall^{\infty}}$ is $\Sigma_5^0$-complete.
\end{thm}

\proof (Sketch.)
The condition
\[ (\forall^{\infty}i)[\TD^+(L^e_i,\cL_e) < \infty] \leftrightarrow
(\exists i)(\forall j \geq i)(\exists a)[\TD^+(L^e_j,\cL_e) < a)], \]
together with the fact that $\TD^+(L^e_j,\cL_e) < a$ is a $\Sigma_3^0$
predicate (as shown in the proof of Theorem \ref{tdplusfinite}),
shows that $I_{TD^+}^{\forall^{\infty}}$ is $\Sigma_5^0$.  The proof that $I_{TD^+}^{\forall^{\infty}}$ is
$\Sigma_5^0$-hard is very similar to that of Theorem \ref{t3complete}.
~\qed

\medskip
\noindent
Finally, consider the positive
extended teaching dimension.  Like the 
u.r.e.\ families with finite
extended teaching dimension, those with finite \emph{positive}
extended teaching dimension have a particularly
simple structure. 
   
\begin{thm}\label{thm:positivextd}
$I_{XTD^+}$ is $\Pi_2^0$-complete.
\end{thm}

\proof (Sketch.)  One may verify directly that
\[ \XTD^+(\cL_e) < \infty \leftrightarrow (\forall i,j)[L_i^e = L_j^e \vee
L^e_i \cap L^e_j = \emptyset].\]

Note that $(\forall i,j)[L_i^e = L_j^e \vee
L^e_i \cap L^e_j = \emptyset] $ iff
\[ (\forall i,j)[(\forall x,s)(\exists t > s)
[\charf_{L^e_{i,t}}(x) = \charf_{L^e_{j,t}}(x)] \vee (\forall s')[L^e_{i,s'} \cap L^e_{j,s'} = \emptyset]],  \]
and that the latter expression reduces to a $\Pi_2^0$ predicate.  Hence $I_{XTD^+}$ is
$\Pi_2^0$.

To establish that $I_{XTD^+}$ is $\Pi_2^0$-hard, take any $\Pi_2^0$ predicate
$P$, and let $g$ be a recursive function such
that $P(e) \leftrightarrow g(e) \in \indset{inf}$.
Let $f$ be a recursive function such that $\cL_{f(e)} = \{\natnum,G\}$,
where
$$
G = \begin{cases}
\natnum & \mbox{if $|W_{g(e)}| = \infty$;} \\
\{0\} \cup \{x: x < m\} & \mbox{if $m$ is the least number}\\
\ & \mbox{such that $(\forall s \geq m)[W_{g(e),s} =$} 
\mbox{$W_{g(e),s+1}]$.}\end{cases}
$$  
Then $P(e) \leftrightarrow g(e) \in \indset{inf} \leftrightarrow \XTD^+(\cL_{f(e)}) < \infty$.~\qed 

\section{Recursive Teaching Dimension}

Although the classical teaching dimension model is quite succinct
and intuitive, it is rather restrictive.  For example, let
$\cL$ be the concept class consisting of the empty set $L_0 = \emptyset$ and
all singleton sets $L_i = \{i\}$ for positive $i\in\natnum$.  
Then $\TD(L_i,\cL) = 1$ for all 
$i \in \natnum\sm\{0\}$ but $\TD(L_0,\cL) = \infty$.  Thus
$\TD(\cL) = \infty$ even though the class $\cL$ is relatively simple.
One deficiency of the teaching dimension model is that it
fails to exploit any property of the learner other than
the learner being \emph{consistent}.  The \emph{recursive teaching
model} \cite{ZillesLHZ11,zeinab}, on the other hand, uses inherent 
structural properties of concept classes to define a
teaching-learning protocol in which the learner is
not just consistent, but also ``cooperates" with the
teacher by learning from a sequence of samples that
is defined relative to the given concept class.  
In this section, we shall consider the arithmetical complexity of the 
index set of the class of all u.r.e.\ families that
are teachable in some variants of the recursive
teaching model.  The complexities of interesting
problems relating to the original recursive teaching
model remain open. 

\begin{defn}\label{defn:rtdoneplus}
A \emph{positive teaching plan for $\mathcal{L}$} is  
is any sequence $\TP = ((L_0,d_0),$ $(L_1,d_1),\linebreak[3]\ldots)$ 
where (i)~the families $\{L_i\}$ form a partition of $\cL$, 
and (ii)~$d_i = \TD^+(L_i,\mathcal{L}\setminus\bigcup_{0\le j<i}\{L_j\})$ for all $i$.
$\RTD^+_1(\mathcal{L})$ is defined to be $\inf\{ord(\TS):\mbox{$\TS$ 
is a positive teaching plan for $\mathcal{L}$}\}$.
Since this paper only considers positive teaching plans, 
positive teaching plans will simply be called ``teaching
plans".  Note that a positive teaching plan for $\cL$ is essentially
a teaching sequence for $\cL$ that employs only positive examples
and partitions $\cL$ into singletons.  
\end{defn}

\noindent We begin with a lemma; the proof is omitted.

\begin{lem}
Let $\{\mathcal{F}_i\}_{i\in\natnum}$ be any sequence of families.
If $\sup\{\rtd{1}{+}{\mathcal F_i} : i \in \bbn\} < \infty$, then $\sup\{\rtd{1}{+}{\mathcal F_i} : i \in \bbn\} \leq \rtd{1}{+}{\bigsqcup_{i\in\bbn} \mathcal F_i} \leq \sup\{\rtd{1}{+}{\mathcal F_i} : i \in \bbn\} + 1$; otherwise, $\rtd{1}{+}{\bigsqcup_{i\in\bbn} \mathcal F_i} = \infty$.
\end{lem}

\begin{defn}
We denote by $R^+_1$ the set of codes for u.r.e.~families, $\mathcal L$, such that $\rtd{1}{+}{\mathcal L}$ is finite.
\end{defn}

\begin{thm}\label{thm:rtdoneplus}
$R^+_1$ is $\Sigma^0_4$-complete.
\end{thm}

\proof To show that $R^+_1$ is $\Sigma^0_4$, fix $e$, a code for a u.r.e.~family, $\mathcal F$.  Given $n \in \bbn$, whether or not $\rtd{1}{+}{\mathcal F} \leq n$ can be decided by executing the following algorithm.  Find the $F \in \mathcal F$ of least index such that $\td{+}{F}{\mathcal F \setminus \{F\}} \leq n$ and call this set $F_0$.  Having defined $F_0, \ldots , F_i$, let $F_{i+1}$ be the set of least index in $\mathcal F \setminus \{F_0, \ldots , F_i\}$ such that $\td{+}{F_{i+1}}{\mathcal F \setminus \{F_0, \ldots , F_{i+1}\}} \leq n$.  If there is a teaching plan for $\mathcal F$ with order at most $n$, then the above algorithm will produce such a teaching plan because 
$\rtd{1}{+}{\mathcal F \setminus \{F_0, \ldots , F_i\}} \leq n$ for all $i \in \bbn$.  
Conversely, if there is no such teaching plan, then clearly the algorithm must initiate a non-terminating search at some stage.

Observe that the statement $\td{+}{F_{i+1}}{\mathcal F \setminus \{F_0, \ldots , F_{i+1}\}} \leq n$ is $\Sigma^0_2$, therefore $\rtd{1}{+}{\mathcal F} \leq n$ is $\Pi^0_3$.  Finally, $\rtd{1}{+}{\mathcal F} < \infty$ is equivalent to $(\exists n) \big($ $\rtd{1}{+}{\mathcal F}\leq n \big)$; hence, it is $\Sigma^0_4$. It remains to show that $R^+_1$ is $\Sigma^0_4$-hard.

Every $\Sigma^0_4$ predicate is of the form $(\forall^{\infty} n) \big( g(e,n) \in \indset{coinf} \big)$, where $g$ is a computable function.  Fix such a predicate, $P$, and computable function $g$.

For $k \in \bbn$, let $f_k: \bbn \rightarrow \bbn$ be a uniformly computable sequence of functions such that 
(1) $D_{f_k(n)} \subsetneq [0,k]$, (2) $(\forall S \subsetneq [0,k])(\exists n) \big( D_{f_k(n)} = S \big)$ 
and (3) $(\forall n \in \mbox{ran}(f_k)) \big( |f_k^{-1}(n)| = \infty \big)$.
Define $\mathcal L_k = \big\{ [0,k] \oplus \emptyset \big\} \cup \big\{ D_{f_k(n)} \oplus \{ n \} : n \in \bbn \big\}$.  We will denote the members of $\mathcal L_k$ by $L^k_i$, where $L^k_0 = [0,k] \oplus \emptyset$ and $L^k_{n+1} = D_{f_k(n)} \oplus \{ n \}$.  Observe that $\rtd{1}{+}{\mathcal L_k} = k+1$. 
However, for any finite $\mathcal L' \subset \mathcal L_k$, $\rtd{1}{+}{\mathcal L'} = 1$.

Let $m$ be a computable function such that $m(a,n,s)$ is the number of $t < s$ such that the $n^{th}$ element of the complement of $W_{a,t}$ differs from the $n^{th}$ element of the complement of $W_{a,t+1}$.  Define
$\mathcal G_{a,i} = \bigsqcup_{n \in \bbn} \big\{  L^i_j : (\exists s)\big( j < m(a,n,s) \big)  \big\}$
and let $\mathcal F_e = \bigsqcup_{i \in \bbn} \mathcal G_{g(e,i),i}$.  By the construction, there is a computable function, $h$, such that $W_{h(e)} = \mathcal F_e$.  We omit the proof that $\rtd{1}{+}{\mathcal F_e} < \infty$ iff $P(e)$.~\qed

\begin{defn}
Given $\sigma \in \bbn^*$ and $S = \{ s_0, \ldots , s_n \} \subset \bbn$ with $s_0 < \cdots < s_n < |\sigma|$, define $\sigma[S] \in \bbn^*$ by $\sigma[S](i) = \sigma(s_i)$ for $i \leq n$.
\end{defn}

\noindent We now consider a ``semi-effective" version of the recursive teaching model
in which the teacher presents only positive teaching sets to the learner. 

\begin{defn} \label{defn:rtdplus}
Let $\mathcal{L}$ be any family of 
subsets of $\mathbb{N}$.  
A \emph{positive teaching sequence for 
$\mathcal{L}$} is any sequence 
$\TS = ((\mathcal{F}_0,d_0),
\linebreak[3](\mathcal{F}_1,d_1),\ldots)$ such that (i)~the 
families $\mathcal{F}_i$ form a partition 
of $\mathcal{L}$ with each $\mathcal{F}_i$ nonempty, 
and (ii) for all $i$ and all
$L \in \mathcal{F}_i$, there is a subset $S_L \subseteq L$
with $|S_L|  = d_i < \infty$ such that for all 
$L' \in \bigcup_{j \geq i}\mathcal{F}_j$,
it holds that $S_L \subseteq L' \rightarrow L = L'$. 
$\sup\{d_i:i\in\mathbb{N}\}$ is called the
\emph{order of $\TS$}, and is denoted by $ord(\TS)$.
The \emph{positive recursive teaching dimension
of $\mathcal{L}$} is defined as $\inf\{ord(\TS):\mbox{$\TS$ is 
a positive teaching sequence for $\mathcal{L}$}\}$ and is denoted by $\RTD^+(\mathcal{L})$.  

We denote by $R^+_{ure}$ the set of codes for u.r.e.~families, $\mathcal L$, such that $\rtd{}{+}{\mathcal L}$ is finite and witnessed by a u.r.e.~teaching sequence.
In this section, a ``teaching sequence" will always mean a
u.r.e.\ teaching sequence. 
\end{defn}

\noindent
Our last major result is that $R^+_{ure}$ is $\Sigma_5^0$-complete, which we establish in the following three theorems.

\begin{thm}\label{rplus-s5-hard}
$R^+_{ure}$ is $\Sigma_5^0$-hard.
\end{thm}

\proof
Fix a $\Sigma_5^0$-predicate $P$.  As in the proof of Theorem 4.4 from \cite{beros-JSL14}, let $g$ be computable such that
$P(e) \rightarrow (\exists x) \Big( (\forall x' > x) (\forall y) \big( g(e,x',y)$ $\in \indset{cof} \big)
\wedge (\forall x' \leq x)(\exists^{\leq 1} y) \big( g(e,x', y) \in \indset{coinf} \big) \Big)$
and $\lneg P(e) \rightarrow (\forall x)(\exists! y) \Big( g(e,x,y)$ $\in \indset{coinf} \Big)$.

As in the proof of Theorem \ref{thm:rtdoneplus}, let $\{f_k\}_{k \in \bbn}$ be a uniformly computable sequence of functions such that for all $k,n \in \bbn$, (1) $D_{f_k(n)} \subsetneq [0,k]$, (2) $(\forall S \subsetneq [0,k])(\exists m) \big( D_{f_k(m)} = S \big)$
and (3) $(\forall m \in \mbox{ran}(f_k)) \big( |f_k^{-1}(m)| = \infty \big)$.

Fix $a,n \in \bbn$ and $\sigma \in \bbn^*$.  Define $H_a(x,\sigma) = \{ \langle \sigma,x \rangle \}$ if $x \in W_a$ and $H_a(x,\sigma) = \emptyset$ otherwise.  Using this notation, we define the set $A^{n} = \bigoplus_{j \in \bbn} ([0,n]$ $\oplus \emptyset)$
and the sets $A^{a,n}_{i,\sigma} = \Big( D_{f_n(\sigma(0))} \oplus H_a(i,\sigma) \Big) \oplus \Big( \bigoplus_{j < i} ([0,n] \oplus \emptyset) \Big)\oplus \Big( \bigoplus_{1 \leq j < |\sigma|}$ $(D_{f_n(\sigma(j))} \oplus \{\langle \sigma[\{0\}\cup[j,|\sigma|)] \rangle\}) \Big)$.
Using the above sets, we define the following families:
$\mathcal G^{a,n}_i = \Big\{ A^{a,n}_{i,\sigma} : \sigma \in \bbn^* \wedge |\sigma| \geq 2 \Big\};
\mathcal G^{a,n} = \Big\{ A^n \Big\} \cup \bigcup_{i \in \bbn} \mathcal G^{a,n}_i$.

Suppose that $a \in \indset{cof}$ and let $x_0, \ldots x_k$ be an increasing enumeration of $\overline{W}_a$.  The following is a teaching sequence for $\mathcal G^{a,n}$:
$\Big( \big( \bigcup_{i \in W_a} \mathcal G^{a,n}_i, 1 \big), \big( \mathcal G^{a,n}_{x_0}, 1 \big), \ldots ,$ 
$\big( \mathcal G^{a,n}_{x_k}, 1 \big), \big( \big\{ A^{a,n} \big\}, 1 \big) \Big).$
Thus, $\rtd{}{+}{\mathcal G^{a,n}} = 1$.  Now suppose that $a \in \indset{coinf}$ and let $x_0, x_1, \ldots$ be an increasing enumeration of $\overline{W}_a$.  Suppose $\TS = ((\mathcal L_0, d_0),$ $(\mathcal L_1, d_1), \ldots )$ is a teaching sequence for $\mathcal G^{a,n}$ and $\order{\TS} \leq n$.  Consider an arbitrary $A^{a,n}_{x_i, \sigma} \in \mathcal G^{a,n}_{x_i}$ for $i \geq 1$.  $A^{a,n}_{x_i, \sigma} \not\in \mathcal L_0$, because $n+1$ points are needed to distinguish $A^{a,n}_{x_i, \sigma}$ from every member of $\mathcal G^{a,n}_{x_0}$.  Since $\mathcal G^{a,n}_{x_1} \cap \mathcal L_0 = \emptyset$, we know that $A^n \not\in \mathcal L_0$.  Now suppose that $\mathcal G^{a,n}_{x_{k}} \subseteq \bigcup_{i \geq k} \mathcal L_i$ and $A^n \not\in \bigcup_{i < k} \mathcal L_i$, then $\mathcal L_{k}$ cannot contain any member of $\mathcal G^{a,n}_{x_{i}}$ for $i \geq k+1$ because $n+1$ points are needed to distinguish the members of $\mathcal G^{a,n}_{x_{k}}$ from the members of $\mathcal G^{a,n}_{x_{i}}$.  As before, this also implies $A^n \not\in \mathcal L_{k+1}$.  By induction, we conclude that $A^n \not\in \mathcal L_i$ for any $i \in \bbn$.  This is a contradiction, so $\rtd{}{+}{\mathcal G^{a,n}} \geq n+1$.  Since
\[
\TS = \Big( \big( \big\{ A^n\big\}, n+1 \big), \big( \mathcal G^{a,n}_0, 1 \big), \big( \mathcal G^{a,n}_1, 1 \big), \ldots \Big)
\]
is a teaching sequence for $\mathcal G^{a,n}$ and $\order{\TS} = n+1$, we conclude that \linebreak[4]$\rtd{}{+}{\mathcal G^{a,n}}$ $= n+1$.  Finally, define
\[
\mathcal F_{e,x} = \bigsqcup_{y \in \bbn} \mathcal G^{g(e,x,y),x} \mbox{\ \ and\ \ } \mathcal F_e = \bigsqcup_{x \in \bbn} \mathcal F_{e,x}.
\]

We wish to prove that $\rtd{}{+}{\mathcal F_e} < \infty$ if and only if $P(e)$.  First, suppose $P(e)$.  For all but finitely many $x$, $g(e,x,y) \in \indset{cof}$ for all $y$.  This means that $\rtd{}{+}{\mathcal F_{e,x}} = 1$ for all but finitely many $x$.  For each $x$ for which $\rtd{}{+}{\mathcal F_{e,x}} \neq 1$ the dimension is still finite, hence, there is a uniform bound, $n$, on the recursive teaching dimension of all the $\mathcal F_{e,x}$.  We conclude that $\rtd{}{+}{\mathcal F_e} < \infty$.

On the other hand, if $\lneg P(e)$, then for every $x$ there is exactly one $y$ such that $g(e,x,y) \in \indset{coinf}$.  Hence, $\mathcal F_e$ is the disjoint union of families whose \RTD\ is unbounded.  We have thus reduced an arbitrary $\Sigma^0_5$-predicate to $R^+_{ure}$.~\qed

\begin{thm}\label{rplus-s5}
$R^+_{ure}$ is $\Sigma_5^0$.
\end{thm}

\proof
Let $\{S_n\}_{n \in \bbn}$ enumerate the u.r.e.~teaching sequences, with $S_n = \big( (\mathcal L^n_0,$ $d^n_0)$, $(\mathcal L^n_1, d^n_1), \ldots \big)$.  Consider a u.r.e.~family, $\mathcal F = \{ F_0, F_1, \ldots \}$ coded by $e$.
\begin{align*}
e \in R^+_{ure} &\leftrightarrow (\exists a,n)\Big( \order{\TS_a} \leq n \wedge\mbox{$\TS_a$ is a teaching sequence for $\mathcal F$}\Big)\\
\order{\TS_a} \leq n &\leftrightarrow (\forall i) \Big( d_i \leq n \Big)
\end{align*}
To say that $S_a$ is a teaching sequence for $\mathcal F$ is equivalent to (1)
$\{\mathcal L_0, \mathcal L_1, \ldots \}$ is a partition of $\mathcal F$ and
(2) $(\forall i\in\bbn, L\in\mathcal L_i)\big( \td{+}{L}{\bigcup_{j \geq i} \mathcal L_j} \leq d_i \big)$.

Since the statement $\td{+}{L}{\bigcup_{j \geq i} \mathcal L_j} \leq d_i$ is $\Sigma^0_3$, we know that the statement $(\forall i\in\bbn, L\in\mathcal L_i)\big( \td{+}{L}{\bigcup_{j \geq i} \mathcal L_j} \leq d_i \big)$ is $\Pi^0_3$.  The statement that $\{\mathcal L_0, \mathcal L_1, \ldots \}$ is a partition of $\mathcal F$ is equivalent to
\begin{align*}
&(\forall i \in \bbn, F \in \mathcal L_i)(\exists F' \in \mathcal F)\Big( F = F' \Big)\wedge (\forall F \in \mathcal F)(\exists i \in \bbn, F' \in \mathcal L_i)\Big( F = F' \Big)\\
\wedge &(\forall i,j \in \bbn, F \in \mathcal L_i, F' \in \mathcal L_j)\Big( i \neq j \rightarrow F \neq F' \Big),
\end{align*}
which is $\Pi^0_4$.
Thus, $e \in R^+_{ure}$ is $\Sigma^0_5$.~\qed

\begin{thm}\label{thm:rtdureplus}
$R^+_{ure}$ is $\Sigma^0_5$-complete.
\end{thm}

\section{Conclusion}

This paper studied the arithmetical complexity of index
sets of classes of u.r.e.\ families that
are teachable under various teaching models.
Our main results are summarised in Table \ref{tab:summary1}.
While u.r.e.\ families constitute a 
very special case of families of sets,
many of our results may be extended to the class of
families of \emph{countably} many sets; more precisely,
if we define $C^X_j = \{x: \spn{j,x} \in X\}$ and
$\cC_X = \{C^X_j: j \in \natnum\}$ for any $X \subseteq \natnum$, 
it is not difficult to apply our results 
to determine the position of $\{X: I(\cC_X) < \infty\}$
for different teaching parameters $I$ in the \emph{hierarchy
of sets of sets} \cite[\S 15.1]{Rog}.
We also determined first-order formulas with the least
possible quantifier complexity defining some
fundamental decision problems in algorithmic teaching.
Our work may be extended in several directions.        
For example, it might be interesting to investigate
the arithmetical complexity of index sets
of classes of general -- even non-u.r.e.\ -- families
that are teachable under the teaching models
considered in the present paper.  In particular,
it may be asked whether the arithmetical complexity
of the class of teachable families with one-one numberings
is less than that of the class of teachable families
that do not have one-one numberings.   

\begin{table}
\begin{scriptsize}
\makebox[\linewidth]{
\renewcommand{\arraystretch}{2}
\begin{tabular}{|c||c|}
\hline
Index Set &\makebox[4cm]{Arithmetical Complexity}\\\hline\hline
$\{e:(\forall L \in \mathcal{L}_e)[\TD^{(+)}(L,\mathcal{L}_e) < \infty]\}$
&$\Pi_4^0$-complete (Thms \ref{t1complete}, \ref{tdpluseachfinite})\\
\hline
$\{e:\TD^{(+)}(\mathcal{L}_e) < \infty\}$
&$\Sigma_5^0$-complete (Thms \ref{thm:td}, \ref{tdplusfinite})\\
\hline
$\{e:\XTD(\mathcal{L}_e) < \infty\}$; $\{e:\XTD^+(\mathcal{L}_e) < \infty\}$
&$\Sigma_3^0$-complete (Thm \ref{thm:xtd}); $\Pi_2^0$-complete (Thm \ref{thm:positivextd})\\
\hline
$\{e:\RTD^+_1(\mathcal{L}_e) < \infty\}$
&$\Sigma_4^0$-complete (Thm \ref{thm:rtdoneplus})\\
\hline
$\{e:\RTD^+_{ure}(\mathcal{L}_e) < \infty\}$
&$\Sigma_5^0$-complete (Thm \ref{thm:rtdureplus})\\
\hline
\end{tabular}
}
\end{scriptsize}
\smallskip
\caption{Summary of main results on u.r.e.\ families. The notation $\TD^{(+)}$ indicates that the result holds for both $\TD$ and $\TD^+$.}\label{tab:summary1}
\end{table}

\end{document}